\documentclass[10pt]{amsart}
\usepackage{amscd}

\theoremstyle{definition}
\newtheorem{thm}{Theorem}[section]
\newtheorem{lem}[thm]{Lemma}
\newtheorem{prp}[thm]{Proposition}
\newtheorem{dfn}[thm]{Definition}
\newtheorem{cor}[thm]{Corollary}

\newtheorem{rmk}[thm]{Remark}

\newcommand{\beq}{\begin{equation}}
\newcommand{\eeq}{\end{equation}}
\newcommand{\beqr}{\begin{eqnarray*}}
\newcommand{\eeqr}{\end{eqnarray*}}
\newcommand{\bal}{\begin{align*}}
\newcommand{\eal}{\end{align*}}
\newcommand{\bei}{\begin{itemize}}
\newcommand{\eei}{\end{itemize}}

\newcommand{\cK}{K}

\newcommand{\af}{\alpha}

\newcommand{\ph}{\varphi}

\newcommand{\Z}{{\mathbb{Z}}}

\newcommand{\C}{{\mathbb{C}}}
\newcommand{\N}{{\mathbb{N}}}

\newcommand{\tsr}{{\mathrm{tsr}}}

\newcommand{\ca}{C*-algebra}

\newcommand{\PSP}{Property~(SP)}

\renewcommand{\S}{\subset}

\newcommand{\Aut}{{\mathrm{Aut}}}

\begin{document}

\title{\Large{Stable rank for inclusions of C*-algebras}}

\author{\Large{Hiroyuki Osaka}}

\date{}

\address{Department of Mathematical Sciences,
 Ritsumeikan University, Kusatsu, Shiga,
 525-8577  Japan} % ???

\email[]{osaka@se.ritsumei.ac.jp}

\thanks{
Research partially supported by the
Open Research Center Project for Private Universities: matching fund
{}from MEXT, 2004-2008, and by a Grant in Aid for Scientific Research,
Ritsumeikan University, 2006.}

%\maketitle

\begin{abstract}
When a unital \ca~$A$ has topological stable rank one 
(write $\tsr(A) = 1$), 
we know that $\tsr(pAp) \leq 1$ for a non-zero projection
$p \in A$. When, however, $\tsr(A) \geq 2$, it is generally faluse. 
We prove that if a unital C*-algebra $A$ has a 
simple unital C*-subalgebra $D$ of $A$ with common unit such that
$D$ has \PSP~ and $\sup_{p\in P(D)}\tsr(pAp) < \infty$, then
$\tsr(A) \leq 2.$ 
As an application let $A$ be a simple unital~ \ca~with $\tsr(A) = 1$ and \PSP, 
$\{G_k\}_{k=1}^n$ finite groups, $\af_k$ actions from $G_k$ to 
${\rm Aut}((\cdots((A\times_{\af_1}G_1)\times_{\af_2} G_2)\cdots)\times_{\af_{k-1}}G_{k-1}).$ 
$(G_0 = \{1\})$
Then 
$$
\tsr((\cdots ((A\times_{\af_1}G_1)\times_{\af_2} G_2)\cdots)\times_{\af_n}G_n)
\leq 2.
$$
\end{abstract}

\maketitle

\section{Introduction}
For a unital \ca~$A,$ recall
that the {\emph{topological stable rank}} $\tsr (A)$ of $A$ is defined
to be the least integer $n$
such that the set ${\mathrm{Lg}}_n (A)$
of all $n$-tuples $(a_1,, a_2 \dots, a_n) \in A^n$
which generate $A$ as a left ideal is dense in $A^n.$
(See Definition~1.4 of~\cite{Rf1}.)
The topological stable rank of a nonunital
\ca\  is defined to be that of its smallest unitization.
Note that $\tsr (A) = 1$ is equivalent to
density of the set of invertible elements in $A.$
Furthermore, $\tsr (A) = 1$ implies that $\tsr (A \otimes M_{n}(\C)) = 1$
for all $n,$ and that $\tsr (A \otimes \cK) = 1,$
where $\cK$ is the algebra of compact
operators on a separable infinite dimensional Hilbert space.
Conversely, if $\tsr (A \otimes M_{n}(\C)) = 1$ for some $n,$
or if $\tsr (A \otimes \cK) = 1,$
then $\tsr (A) = 1.$
(See Theorems 3.3 and~3.6 of~\cite{Rf1}.)
Simple AH~algebras with slow dimension growth have
topological stable rank one (Theorem~1 of~\cite{BDR}),
as do irrational rotation algebras (\cite{pu}).
If $A$ is unital and $\tsr (A) = 1,$
then $A$ has cancellation (Proposition 6.5.1 of~\cite{Bl1}).
It follows immediately that $A$ is stably finite in the sense
that no matrix algebra $M_n (A)$ has an infinite projection.

Blackadar proposed the question in \cite{bl3} whether 
$\tsr(A \times_\af G) = 1$ for any unital AF \ca\ $A$, 
a finite group $G$, and an action $\af$ of $G$ on $A.$
Osaka and Teruya \cite{OT2} gave a partial answer to a generalized 
Blackadar's question, that is, if $A$ is a simple unital 
\ca\ with $\tsr(A) = 1$ and \PSP\, then 
$\tsr(A \times_\af G) \leq 2$ for a finite group $G$, 
and an action $\af$ of $G$ on $A.$ 
We do not know whether $\tsr(A \times_\af G) = 1$ then.
He also studied about topological stable rank of corner algebra 
in \cite{bl4} and showed that if  a non-zero projection $p$ is 
full in a given C*-algebra $A$, then 
$\tsr(A) \leq \tsr(pAp).$ When $\tsr(A) = 1$, $\tsr(pAp) = 1$ for 
any non-zero projection by \cite{Rf1}. 

In this paper we show that the boundedness of a set 
 $\{\tsr(pAp)\colon p\ \hbox{projection}\ \in D\}$ implies that 
$\tsr(A) \leq 2$ when a unital \ca\ $A$ has a unital simple 
C*-subalgebra $D$ of $A$ with $\tsr(D) = 1$ and 
\PSP. 
As an application  let $A$ be a simple unital~ \ca~with $\tsr(A) = 1$ and \PSP, 
$\{G_k\}_{k=1}^n$ finite groups, and $\af_k$ actions from $G_k$ to 
${\rm Aut}((\cdots((A\times_{\af_1}G_1)\times_{\af_2} G_2)\cdots)\times_{\af_{k-1}}G_{k-1}).$ 
$(G_0 = \{1\})$
Then 
$$
\tsr((\cdots ((A\times_{\af_1}G_1)\times_{\af_2} G_2)\cdots)\times_{\af_n}G_n)
\leq 2.
$$

The author would like to thank Tamotsu Teruya for a fruitful discussion.

\section{Stable rank of corner algebras}

Let $P(A)$ be a set of all projections in a unital C*-algebra $A.$
When a unital C*-algebra $A$ has topological stable rank one, we have $\tsr(pAp) = 1$ 
for any non-zero projection $p$ in $A$. But when $A$ has topological stable rank 
more than one, we do not know whether the set of $\{\tsr(pAp)\colon p \in p(A)\}$
is bounded. The following implies that the boundedness of a set 
$\{\tsr(pAp)\colon p \in p(A)\}$ controls topological stable rank of a given C*-algebra.

Recall that a \ca~$A$ has 
{\emph{\PSP}} if there exists a non-zero projection
in any non-zero hereditary subalgebra of $A$.

\begin{prp}\label{P:Stablerank2}
Let $1 \in D \subset A$ be an inclusion of 
\ca s with common unit. Suppose that
\begin{align*}
&(i) D \ \hbox{is simple}\\
&(ii) D \  \hbox{has \PSP}\\
&(iii) \sup_{p\in P(D)}\tsr(pAp) < \infty
\end{align*}
Then $\tsr(A) \leq 2.$
\end{prp}

\begin{proof}
Set $K = \sup_{p\in P(D)}\tsr(pAp).$   
Since $D$ is simple with \PSP,
there exist $N \in \N$ such that $N > K$ and mutually orthogonal projections 
$\{e_j\}_{j=1}^N \subset D$
such that $e_i \sim e_j$ by Lemma~ 3.5.7 of \cite{Hl}.

Set $e = \sum_{j=1}^Ne_j$. 
Then $eAe \cong M_N(e_1Ae_1).$
Hence from Theorem~6.1 of \cite{Rf1}
\begin{align*}
\tsr(eAe) &= \tsr(M_N(e_1Ae_1))\\
&= \{\dfrac{\tsr(e_1Ae_1) - 1}{N}\} + 1\\
&\leq 2,
\end{align*}
where $\{t \}$ denotes the least integer which is greater than
or equal to $t$. 
Since $e$ is full in $D$, so is $e$ in $A.$ 
From Theorem~4.5 of \cite{bl4} 
\begin{align*}
\tsr(A) \leq \tsr(eAe) \leq 2.
\end{align*}
\end{proof}

\vskip 3mm

\begin{dfn}\label{D:FiniteDepth}
Let $1 \in A \subset B$ be an inclusion of unital \ca s
with a conditional expectation $E \colon B \to A$
with index-finite type.
That is, there exists a set $\{(u_i, u_i^*)\}_{i=1}^n
\subset B \times B$ such that 
$$
x = \sum_{i=1}^nu_iE(u_i^*x) = \sum_{i=1}^nE(xu_i)u_i^*
$$
for any $x \in B$. We call a set $\{(u_i, u_i^*)\}_{i=1}^n$ 
{\emph{a quasi-basis}} for $E$ (Definition 1.~2.~2 of \cite{wata}).
Set $B_0 = A,$ $B_1 = B,$
and $E_1 = E.$
Recall the \ca\  version of the basic construction
(Definition~2.2.10 of~\cite{wata}, where it is called the
C*~basic construction).
We inductively define $e_k = e_{B_{k - 1}}$
and $B_{k + 1} = C^* (B_k, e_k),$
the Jones projection and \ca\  for the basic construction
applied to $E_k \colon B_k \to B_{k - 1},$
and take $E_{k + 1} \colon B_{k + 1} \to B_k$ to be
the dual conditional expectation $E_{B_k}$ of
Definition~2.3.3 of~\cite{wata}.
This gives the {\emph{tower of iterated basic constructions}}
\[
B_0 \subset B_1 \subset B_2 \subset \cdots \subset B_k \subset \cdots,
\]
with $B_0 = A$ and $B_1 = B.$
It follows from Proposition 2.10.11 of~\cite{wata} that
this tower does not depend on the choice of $E.$

We then say that the inclusion $A \subset B$
has {\emph{finite depth}}
if there is
$n \in \N$ such that
$(A' \cap B_{n}) e_{n} (A' \cap B_{n}) = A' \cap B_{n + 1}.$
We call the least such $n$ the {\emph{depth}} of the inclusion.
\end{dfn}

\vskip 3mm

\begin{lem}\label{Stablerankbounded}
Let $B$ be a unital \ca\ and 
let $D \S A \S B$ be unital subalgebras with common unit.
Let $E \colon B \to A$ be
a faithful conditional expectation with index-finite type
and depth 2, that is, 
$(A' \cap B_2)e_2(A' \cap B_2) = A' \cap B_3.$
Suppose that $\sup_{p\in P(D)}\tsr(pAp) < \infty.$ 
Then $\sup_{p\in P(D)}\tsr(pBp) < \infty.$
\end{lem}

\begin{proof}
Set $K = \sup_{p\in P(D)}\tsr(pAp).$
Since an inclusion $A \S B$ is of index finite and depth 2, 
there exist $n \in \N$ and  a quasi-basis $\{(v_k,v_k^*)\}_{k=1}^n \in 
(A' \cap B_2) \times (A' \cap B_2)$ for $E_2$ by 
Theorem 4.4 of \cite{OT}. 

Let $A \S B \S B_2$ be a basic construction, 
$\{(u_j,u_j^*)\}_{k=1}^m$ be a quasi-basis for $E$, 
and 
$\ph\colon B_2 \rightarrow qM_m(A)q$ be a canonical 
isomorphism in the argument before Lemma~3.3.4 of \cite{wata} by 
$\ph(xe_Ay) = [E(u_i^*x)E(yu_j)]$, where 
$q = [E(u_i^*u_j)]$ is a projection in $M_m(A)$.

Let $p \in D$ be non-zero projection and 
set $r = [E(u_i^*)E(u_j)].$
Note that since each $v_k$ $(1 \leq k \leq n)$ belongs to 
$A' \cap B_2$,
the restricted map $F$ of $E_2$ to $pB_2p$ is a 
faithful conditional expectation from $pB_2p$ to 
$pBp$ and $\{(pv_k,pv_k^*)\}_{k=1}^n$ is a quasi-basis 
for $F$. 
Then we have 
\begin{align*}
rq &= [E(u_i^*)E(u_j)][E(u_i^*u_j)]\\
&= [\sum_kE(u_i^*)E(u_k)E(u_k^*u_j)]\\
&= [E(u_i^*)E(u_j)] = r
\end{align*}
and
%$$
%\begin{array}{ll}
\begin{align*}
\ph(p) &= \ph(\sum_kpu_ke_Au_k^*)\\
&= \sum_k\ph(pu_ke_Au_k^*)\\
&= \sum_k[E(u_i^*pu_k)E(u_k^*u_j)]\\
&= [E(u_i^*p\sum_ku_kE(u_k^*u_j))]\\
&= [E(u_i^*pu_j)]\\
\end{align*}
% \end{array}
%$$
Since  $r$ is full in $\ph(B_2)$ and 

%$$
%\begin{array}{ll}
\begin{align*}
r\ph(p) &= [E(u_i^*)E(u_j)][E(u_i^*pu_j)]\\
&= [\sum_kE(u_i^*)E(u_k)E(u_k^*pu_j)]\\
&= [E(u_i^*)E(\sum_ku_kE(u_k^*pu_j))]\\
&= [E(u_i^*)E(pu_j)]\\
&= [E(u_i^*)pE(u_j)]\\
&= \ph(p)r
\end{align*}
%\end{array}
%$$
we have 
\begin{align*}
\tsr(pBp) &= \tsr(\ph(p)\ph(B)\ph(p))\\
&\leq \tsr(\ph(p)\ph(B_2)\ph(p)) + n - 1\\
&\leq \tsr(r\ph(p)\ph(B_2)\ph(p)r) + n - 1\\
\end{align*}
The second inequality comes from 
Theorem 2.2 of \cite{JOPT} and 
the third inequality comes from 
Theorem~4.5 of \cite{bl4}.

Set $X = \left(\begin{array}{cccc}
E(u_1)&E(u_2)&\cdots&E(u_m)\\
0&0&\cdots&0\\
\vdots&\vdots&&\vdots\\
0&0&\cdots&0
\end{array}
\right)$. Then $X^*X = r$ and 
a map ${\rm Ad}(X)\colon r\ph(p)\ph(B_2)\ph(p)r \rightarrow pAp$
induces an isomorphism. Indeed, from the same argument in Cuntz \cite{Cu}
\begin{align*}
{\rm Ad}(X)(r\ph(p)\ph(B_2)\ph(p)r) &= Xr\ph(p)\ph(B_2)\ph(p)rX^*\\  
&= (XX^*)X\ph(p)\ph(B_2)(XX^*)X\ph(p)X^*(XX^*)\\
&= 
\left(
\begin{array}{cccc}
1&0&\cdots&0\\
0&0&\cdots&0\\
\vdots&\vdots&&\vdots\\
0&0&\cdots&0
\end{array}
\right)
\left(
\begin{array}{cccc}
p&0&\cdots&0\\
0&0&\cdots&0\\
\vdots&\vdots&&\vdots\\
0&0&\cdots&0
\end{array}
\right)
X\ph(B_2)X^*\\
&\times 
\left(
\begin{array}{cccc}
p&0&\cdots&0\\
0&0&\cdots&0\\
\vdots&\vdots&&\vdots\\
0&0&\cdots&0
\end{array}
\right)
\left(
\begin{array}{cccc}
1&0&\cdots&0\\
0&0&\cdots&0\\
\vdots&\vdots&&\vdots\\
0&0&\cdots&0
\end{array}
\right)
\\
&= 
\left(
\begin{array}{cccc}
p&0&\cdots&0\\
0&0&\cdots&0\\
\vdots&\vdots&&\vdots\\
0&0&\cdots&0
\end{array}
\right)
XX^*\ph(B_2)XX^*
\left(
\begin{array}{cccc}
p&0&\cdots&0\\
0&0&\cdots&0\\
\vdots&\vdots&&\vdots\\
0&0&\cdots&0
\end{array}
\right)\\
&=
\left(
\begin{array}{cccc}
p&0&\cdots&0\\
0&0&\cdots&0\\
\vdots&\vdots&&\vdots\\
0&0&\cdots&0
\end{array}
\right)
\left(
\begin{array}{cccc}
A&0&\cdots&0\\
0&0&\cdots&0\\
\vdots&\vdots&&\vdots\\
0&0&\cdots&0
\end{array}
\right)
\left(
\begin{array}{cccc}
p&0&\cdots&0\\
0&0&\cdots&0\\
\vdots&\vdots&&\vdots\\
0&0&\cdots&0
\end{array}
\right)\\
&= \left(
\begin{array}{cccc}
pAp&0&\cdots&0\\
0&0&\cdots&0\\
\vdots&\vdots&&\vdots\\
0&0&\cdots&0
\end{array}
\right)\\
&\cong pAp.
\end{align*}

Therefore we have
\begin{align*}
\tsr(pBp) &\leq \tsr(pAp) + n - 1\\
&\leq K + n - 1.
\end{align*}
\end{proof}

%\vskip 5mm
\newpage

\section{Main theorem}

\begin{thm}\label{T:Stablerank2}
Let $B$ be a unital \ca,
let $D \S A \S B$ be unital subalgebras,
let $E \colon B \to A$ be
a faithful conditional expectation with index-finite type
and depth 2. 
Suppose that \begin{align*}
&(i) D \ \hbox{is simple}\\
&(ii) D \  \hbox{has \PSP}\\
&(iii) \sup_{p\in P(D)}\tsr(pAp) < \infty
\end{align*}
Then $\tsr(B) \leq 2.$
\end{thm}

\begin{proof}
From Lemma~\ref{Stablerankbounded} 
$\sup_{p\in P(D)}\tsr(pBp) < \infty.$ 
Therefore since $D$ is a unital C*-subalgebra of $B$ 
satisfying three conditions in Proposition~\ref{P:Stablerank2},
we get the conclusion.
\end{proof}

\vskip 5mm

\begin{cor}\label{Crossedproducts}
Let $A$ be a simple unital~ \ca~with $\tsr(A) = 1$ and \PSP, 
$\{G_k\}_{k=1}^n$ finite groups, $\af_k$ actions from $G_k$ to 
${\rm Aut}((\cdots((A\times_{\af_1}G_1)\times_{\af_2} G_2)\cdots)\times_{\af_{k-1}}G_{k-1}).$ 
$(G_0 = \{1\})$
Then 
$$
\tsr((\cdots ((A\times_{\af_1}G_1)\times_{\af_2} G_2)\cdots)\times_{\af_n}G_n)
\leq 2.
$$
\end{cor}

\begin{proof}
Note that since $A$ is simple with $\tsr(A) = 1$, 
we have $\sup_{p\in A}\tsr(p(A\times_{\af_1}G_1)p) < \infty$
from the same argument as in the proof of Theorem 5.1 of \cite{OT}.

Since $1 \in A \subset A\times_{\af_1}G_1 \subset 
(A \times_{\af_1}G_1)\times_{\af_2}G_2$  of finite index and depth 2 
from Lemma 3.1 of \cite{OT}.
Therefore, from Theorem~\ref{T:Stablerank2} we conclude that 
$$
\tsr((A \times_{\af_1}G_1)\times_{\af_2}G_2) \leq 2.
$$
Note that 
$$
\sup_{p\in A}\tsr(p((A \times_{\af_1}G_1)\times_{\af_2}G_2)p) < \infty
$$
from the argument in the proof of Lemma~\ref{Stablerankbounded}.

By induction steps we get the conclusion.
\end{proof}

\vskip 5mm

\begin{rmk}\label{blackadarquestion}
In \cite{bl3} Blackadar presented a question as follows:
Let $A$ be an AF algebra, $G$ a finite group, and $\af$ 
an action of $G$ on $A$. Is it true that $\tsr(A \times_\af G) = 1 ?$
This is still an open question. 
Related to this question we may consider a general version of this. 
That is, let $A$ be a simple \ca\ with $\tsr(A) = 1$, 
$G$ a finite group, and $\af$ 
an action of $G$ on $A$. Is it true that $\tsr(A \times_\af G) = 1 ?$
In this version we can not drop simplicity of $A$. 
See an example in Remark~\ref{blackadarexample} below.  

Corollary~\ref{Crossedproducts} implies that topological stable ranks of 
C*-algebras constructed through crossed product ways by finite groups 
of a simple C*-algebra with topological stable rank one and 
property \PSP\ are less than or equal to 2. 
Note that in this case a crossed product has {\emph{Cancellation Property}} 
by Theorem~4.6 of \cite{JOPT}.
\end{rmk}

\begin{dfn}
Let $A$ be a unital simple AH algebra. 
We say that $A$ has slow dimension growth if 
$A = \lim_{n\rightarrow \infty}A_n$, where 
$A_n = \oplus_{j=1}^{i(n)}P_{(n,j)}M_{n(j)}(C(X_{(n,j)}))P_{(n,j)}$,
$X_{(n,j)}$ is a connected finite CW complex, 
$P_{(n,j)}$ is a projection in  $M_{n(j)}(C(X_{(n,j)}))$, and 
\begin{align*}
\lim_{n\rightarrow\infty}\max_j
\{\frac{\dim X_{(n,j)}}{{\rm rank}P_{(n,j)}}\} = 0.
\end{align*}
\end{dfn}

\vskip 5mm

\begin{cor}\label{TsrAHalgebras}
Let $A$ be a simple, unital, AH-algebra of real rank zero 
which has slow dimension growth and  
$B$ be a unital \ca~with $\tsr(B) < \infty$. 
Let 
$\{G_k\}_{k=1}^n$ be finite groups and $\af_k$ actions from $G_k$ to 
${\rm Aut}((\cdots((B \otimes A) \times_{\af_1}G_1)\times_{\af_2}G_2)\cdots)\times_{\af_{k-1}}G_{k-1}).$ 
$(G_0 = \{1\})$
Then 
$$
\tsr((\cdots ((B \otimes A)\times_{\af_1}G_1)\times_{\af_2} G_2)\cdots)\times_{\af_k}G_n)
\leq 2.
$$
In particular, 
$$
\tsr(B \otimes A) \leq 2.
$$
\end{cor}

\begin{proof}
Let $X$ be a compact Hausdorff space. Then 
$$
\tsr(B \otimes C(X)) \leq \tsr(B) + \dim X
$$
from \cite[Corollary~7.2]{Rf1} or \cite[Theorem~1.13]{nop}.
Hence if $X$ is a connected compact Hausdorff space and 
$p$ is a projection in $M_n(C(X))$, then we have 
\begin{align*}
\tsr(B \otimes pM_n(C(X))p) 
&= \tsr(B \otimes M_{{\rm rank}(p)}(C(X)))\\
&= \tsr(M_{{\rm rank}(p)}(B \otimes C(X)))\\
&= \{\frac{\tsr(B \otimes C(X)) - 1}{{\rm rank}(p)}\} + 1\\
&\leq \frac{\tsr(B) + \dim X - 1}{{\rm rank}(p)} + 1\\
&\leq \frac{\dim X}{{\rm rank}(p)} + \tsr(B) - 1 + 1\\
&\leq \tsr(B) + 1.
\end{align*}
if $\dim X < {\rm rank}(p)$. 
The third equality comes from  
Theorem~6.1 of \cite{Rf1}.

Therefore, since $A$ has slow dimension growth, we have 
$$
\sup_{p \in P(A)}\tsr(B \otimes pAp) \leq \tsr(B) + 1
$$

Let $D = 1 \otimes A$. 
Since $A$ has real rank zero, $A$ has \PSP\ by 
\cite[Theorem 2.6]{bp0}. Hence 
$D \subset B \otimes A$ and $D$ satisfies three conditions in 
Theorem~\ref{T:Stablerank2}. 
Hence as in the same argument in the proof of 
Corollary~\ref{Crossedproducts} we get the conclusion.

\end{proof}

\vskip 5mm

\begin{rmk}\label{blackadarexample}
The estimate in Corollary~\ref{TsrAHalgebras} is 
best possible. 
Indeed, in Example~8.2.1 of \cite{bl3} 
there exists a symmetry $\af$ on 
$C[0, 1] \otimes UHF$ such that 
$\tsr((C[0, 1] \otimes UHF) \times_\af\Z/2\Z) = 2$.
\end{rmk}

\vskip 5mm

\begin{rmk}
Let $A$ be a unital C*-algebra with $\tsr(A) \geq 2$.
Then we can not hope that 
$\tsr(pAp) \leq 2$ for any non-zero projection $p \in A$. 
For example, consider $A = M_n(C[0,1]^n).$ Then 
$\tsr(A) = \{\frac{n-1}{n}\} + 1 = 2$ by Theorem~6.1 of 
\cite{Rf1}. Take $p = {\rm diag}(1,0,\cdots,0) \in M_n(A).$ Then
$p$ is a projection and $\tsr(pAp) = \tsr(C[0,1]^n) = [\frac{n}{2}] + 1$ 
by \cite{LNV}.  
If we take $n \geq 4$, then $\tsr(pAp) \geq 3$.
\end{rmk}

\vskip 5mm

\begin{rmk}
Let $A$ be a unital C*-algebra and $\alpha$ be an action from 
a finite group $G$ to $\Aut(A)$ with $\tsr(A) \geq 2$. 
 We can not hope that $\tsr(A \times_\alpha G) \geq 2$. 
 Indeed, when $A = (C[0, 1] \otimes UHF) \times_\alpha\Z/2\Z$, 
 where $\alpha$ is a symmetry in  Remark \ref{blackadarexample}. Then 
$\tsr(A) = 2$. Let $\hat{\alpha}$ is a dual action of $\alpha$. 
From Takai's duality Theorem (\cite{Ta}) we have 
$$A \times_{\hat{\alpha}}\Z/2\Z 
\cong (C[0, 1] \otimes UHF) \otimes M_2(\C).
$$
Therefore, $\tsr(A \times_{\hat{\alpha}}\Z/2\Z) = 1$. 

On the contrary, 
when $\tsr(A) \ge 2n$, we could  conclude that 
$\tsr(A \times_\alpha \Z/2\Z) \geq n$. 

\begin{prp}
Let $A$ be a unital C*-algebra with $\tsr(A) \geq 2n$. 
Then for any symmetry $\alpha$ on $A$ 
we have $\tsr(A \times_\alpha \Z/2\Z) \geq n$.
\end{prp}

\begin{proof}
From \cite[Theorem~7.1]{Rf1} we have 
$$
\tsr((A \times_\alpha  \Z/2\Z)\times_{\hat{\alpha}} \Z/2\Z)
\leq \tsr(A \times_\alpha  \Z/2\Z) + 1.
$$

On the contrary, from Takai's duality Theorem (\cite{Ta}) 
and \cite[Theorem~6.1]{Rf1} 
\begin{align*}
\tsr((A \times_\alpha  \Z/2\Z)\times_{\hat{\alpha}} \Z/2\Z) 
&= \tsr(A \otimes M_2(\C))\\
&= \{\frac{\tsr(A) - 1}{2}\} + 1\\
&\geq \{\frac{2n-1}{2}\} + 1 = n + 1.
\end{align*}

Therefore, we have 
$$
\tsr(A \times_\alpha  \Z/2\Z) + 1 \geq n + 1,
$$
that is,
$$
\tsr(A \times_\alpha  \Z/2\Z) \geq n.
$$
\end{proof}
\end{rmk}

\end{document}